\documentclass[10pt,letter]{article}

\usepackage{lipsum}

\usepackage{tikz-cd}
\usepackage{amsmath}
\usepackage{amssymb}
\usepackage{amsthm, thmtools}
\usepackage{mathtools}
\usepackage{mathrsfs}

\usepackage{bbm}
\usepackage{bm}
\usepackage{MnSymbol}   
\usepackage{cases}
\usepackage[normalem]{ulem}

\usepackage[width=.80\textwidth,font=footnotesize]{caption}
\usepackage[letterpaper, margin=1.25in]{geometry}
\usepackage[title]{appendix}
\usepackage{fancyhdr}
\usepackage{caption}
\usepackage{subcaption}

\usepackage[shortlabels]{enumitem}
\usepackage[ddmmyyyy]{datetime}
\usepackage{hhline}

\usepackage{listings}
\newcommand{\N}{\mathbb{N}}
\newcommand{\Z}{\mathbb{Z}}
\newcommand{\R}{\mathbb{R}}

\DeclarePairedDelimiter\abs{|}{|}
\DeclarePairedDelimiterX{\inner}[2]{\langle}{\rangle}{#1,#2}
\DeclarePairedDelimiterX{\norm}[1]{\|}{\|}{#1}

\DeclareMathOperator{\diam}{diam}

\DeclareMathOperator{\ExpOp}{\mathbb{E}}
\newcommand{\E}[1]{\ExpOp\left[#1\right]}

\DeclarePairedDelimiter\floor{\lfloor}{\rfloor}

\newtheorem{theorem}{Theorem}[section]
\newtheorem*{theorem*}{Theorem}
\newtheorem{lemma}[theorem]{Lemma}	\newtheorem*{lemma*}{Lemma}
\newtheorem*{cor*}{Corollary}

\theoremstyle{definition}	
\newtheorem{definition}[theorem]{Definition}
\newtheorem{assumption}{Assumption}

\theoremstyle{definition}

\theoremstyle{plain}		
\theoremstyle{definition}	\newtheorem{remark}[theorem]{Remark}

\theoremstyle{definition}	

\theoremstyle{definition}

\begin{document}
\begin{center}
{\Large Existence, Regularity, and a Strong It{\^o}~Formula for the Isochronal Phase of SPDEs} \\
{\textsc{ (preprint)} }
\end{center}

\begin{center}
\begin{minipage}[t]{.75\textwidth}
\begin{center}
Zachary P.~Adams \\
\footnotesize{Freie Universit{\"a}t Berlin, Scads.AI Leipzig,\\ \&~Max Planck Institute for Mathematics in the Sciences}\\
zachary.adams@mis.mpg.de\\
$\,$\\
\hfill
\end{center}
\end{minipage}
\end{center}
\begin{center}
\today 
\end{center}
\vspace{.5cm}
\setlength{\unitlength}{1in}

\vskip.15in

{
\centerline
{\large \bfseries \scshape Abstract}
We prove the existence and regularity of the isochron map for stable invariant manifolds of a large class of evolution equations. 
Our results apply in particular to the isochron map of reaction-diffusion equations and neural field equations. 
Using the regularity properties proven here, we prove a strong It{\^o}~formula for the isochronal phase of stochastically perturbed travelling waves and other patterns appearing in SPDEs driven by white noise, even for SPDEs that only admit mild solutions. 
}\\

\hfill

\section{Introduction}
\label{sec:Intro}
This note presents simple arguments on the existence and regularity of the isochron map for a class of dynamical systems $(\phi_t)_{t\ge0}$ on a separable Banach space in the vicinity of a stable, normally hyperbolic invariant manifold $\Gamma$. 
The isochron map was introduced for stable limit cycles of dynamical systems in $\R^d$ by \cite{G75}~and \cite{W74}. 
Our framework encompasses infinite dimensional dynamical systems, in particular reaction-diffusion systems and neural field equations. 
From the regularity properties of the isochron map proven here, we obtain a strong It{\^o}~formula for the isochronal phase of stochastic perturbations of reaction-diffusion systems, even when these systems only admit a mild solution theory. 

The proofs in this note are short and general, and recover results on the isochron map proven elsewhere in the literature using different techniques \cite{AM22}. 
As will be seen, while \cite{AM22}~requires the invariant manifold to be smooth, our arguments only require $\Gamma$ to be $C^1$. 
Additionally, we make no assumptions on the spectral properties of the generator of the linearization of $\phi_t$ about $\Gamma$. 
This is of particular interest in the case where $\Gamma$ is the orbit of a non-commutative Lie group (such as $SE(n)$ when $\Gamma$ consists of spiral wave solutions to a parabolic PDE), where trying to track the linearized dynamics along $\Gamma$ can lead to significant complications, as noted in the introduction of \cite{KMZ22}. 

As a consequence, the approach taken in this note may lead to significant generalizations of results obtained for travelling pulses and other transient patterns appearing in stochastic evolution equations using other notions of ``phase''. 
For instance, those studied in neural field equations by  \cite{IM16,L16,LS16}, and in reaction-diffusion equations by \cite{EGK20,HH19, HH20a, HH20b,KS14,KS17,S13}. 
The approach of this note also allows us to prove a strong It{\^o}~formula for the isochronal phase of SPDEs with only mild solutions, 
which is stronger than the mild It{\^o}~formula proven using more complicated techniques in \cite{DPJR19,I84}. 

\section{Setup} 
\label{sec:Setup}
We consider an evolution equation of the form 
\begin{equation}
\label{eq:ODE}
\partial_tx\,=\, Lx + N(x)\,\eqqcolon\,V(x), 
\end{equation}
and its stochastic perturbation 
\begin{equation}
\label{eq:SDE}
dX\,=\, \left(LX + N(X)\right)\,dt + \sigma B\,dW. 
\end{equation}
Here, $L$ is a linear operator generating a $C_0$-semigroup on a Hilbert space $H$, $N$ is a nonlinearity defined on a subspace of $H$, $B$ is a linear operator on $H$, $(W_t)_{t\ge0}$ is a cylindrical Wiener process, as in \cite[Chapter 4.1.2]{DPZ14}, and $\sigma\ge0$. 
We make the following assumptions on the deterministic dynamics of \eqref{eq:ODE}. 

\begin{assumption}
\label{assn:PDE}
The system \eqref{eq:ODE}~satisfies the following. 
\begin{enumerate}[(a)]
\item $N$ is a nonlinearity defined on a dense subset $E$ of $H$ such that $E$ is a Banach space with norm $\norm*{\cdot}_E$, and the embedding $E\hookrightarrow H$ is dense and continuous. 
$N$ is 
four times continuously Fr{\'e}chet differentiable 
in the topology of $E$, with first and second Fr{\'e}chet derivatives in this topology denoted $DN$ and $D^2N$. 
\item $L$ is a linear operator with a dense domain of definition $D(L)$ in $H$, and generates a $C_0$-semigroup $(\Lambda_t)_{t\ge0}$ on $H$ that restricts to $E$. 
Moreover, letting $\norm*{\cdot}$ denote the norm of $H$, there exists $\omega>0$ such that $\sup\left\{\norm*{\Lambda_t},\, \norm*{\Lambda_t}_E\right\}\,\le\,e^{-\omega t}$ for all $t\ge0$. 
\end{enumerate}
\end{assumption}

Let the flow map of \eqref{eq:ODE}~be $(t,x)\,\mapsto\,\phi_t(x)$, and suppose that the flow can be expressed by the variation of constants formula as 
\[
\phi_t(x)\,=\, \Lambda_t x + \int_0^t \Lambda_{t-s}N(\phi_s(x))\,ds. 
\]
Since the nonlinearity $N$ is $C^4$ in the topology of $E$, the flow map is $C^3$ in this topology. 
The first and second 
derivatives of $\phi_t$ at $x_0$ in directions $x,y\in E$ are denoted 
\begin{equation}
\label{eq:Derivatives}
x\,\mapsto\,D\phi_t(x_0)[x],\qquad (x,y)\,\mapsto\,D^2\phi_t(x_0)[x,y]. 
\end{equation}

We now state our assumptions on the flow of \eqref{eq:ODE}. 
Specifically, we assume the existence of a stable invariant manifold $\Gamma$ of \eqref{eq:ODE}~that models spatiotemporal patterns. 
In particular, the following guarantees that the dynamics of \eqref{eq:ODE}~is not chaotic on $\Gamma$.

\begin{assumption}
\label{assn:Manifold}
The deterministic system \eqref{eq:ODE}~has a stable, finite dimensional, normally hyperbolic invariant manifold $\Gamma$, as defined in \cite[Condition (H3)]{BLZ98}. 
$\Gamma$ is parameterized by a manifold $\mathcal{S}\subset\R^m$, where $m\in\N$. 
We write $\Gamma=\{\gamma_\alpha\}_{\alpha\in\mathcal{S}}$, and let $B(\Gamma)$ denote the basin of attraction of $\Gamma$ in $H$. 
Moreover, the following hold. 
\begin{enumerate}[(a)]
\item $\phi_t$ and $\phi_t^{-1}$ are Lipschitz on $\Gamma$ with uniform in time Lipschitz constants. 
\item $D\phi_t(\gamma_\alpha)$ is invertible with uniformly bounded inverse for all $\alpha\in\mathcal{S}$. 
\item The map $\alpha\mapsto\gamma_\alpha$ is continuously differentiable with derivative $D\gamma_\alpha$, and $D\gamma_\alpha$ is invertible with a uniformly bounded inverse. 
\end{enumerate}
\end{assumption}

\begin{definition}
\label{def:Isochron}
Assumption \ref{assn:Manifold}~implies that for each $x\in B(\Gamma)$ there is a unique $\pi(x)\in\mathcal{S}$ such that 
\begin{equation}
\label{eq:pidefinition}
\norm*{\phi_t(x)-\phi_t(\gamma_{\pi(x)})}\,\xrightarrow[t\rightarrow\infty]{}\,0, 
\end{equation} 
providing a well-defined map $\pi:B(\Gamma)\rightarrow\mathcal{S}$. 
The existence and uniqueness of this map is proven in Theorem \ref{theorem:C2}. 
We refer to $\pi$ as the \emph{isochron map}~of $\Gamma$. 
\end{definition}

Throughout this note, for $\delta>0$ we define $\Gamma_\delta\,\coloneqq\,\left\{x\in E\,:\,\norm*{x-\gamma_{\pi(x)}}_E\le\delta\right\}\subset E$. 
%
When considering the SDE \eqref{eq:SDE}, we require the following additional assumptions. 

\begin{assumption}
\label{assn:Noise}
Assumption \ref{assn:PDE}~holds, and the SDE \eqref{eq:SDE}~satisfies the following. 
\begin{enumerate}[(a)] 
\item There exists a unique $E$-valued mild solution $(X_t)_{t\ge0}$ to \eqref{eq:SDE}, satisfying 
\begin{equation}
\label{eq:Zmild}
X_t\,=\,\Lambda_tX_0+\int_0^t\Lambda_{t-s}N(X_s)\,ds + \sigma\int_0^t\Lambda_{t-s}B\,dW_s,\qquad X_0\in H, 
\end{equation}
for $X_0\in E$ and all $t<\tau_{blowup}$, where $\tau_{blowup}$ is in $(0,\infty]$ almost surely. 
\item 
For $x\in H$, $(\phi_t(x))_{t>0}$ is continuous in $D(L)$, equipped with the graph norm of $L$. 
\item 
There exists an orthonormal basis $\{e_k\}_{k\in\N}$ of $H$ consisting of eigenfunctions of $L$, and $e_k\in E$ for $k\in\N$. 
Letting $\lambda_k$ denote the eigenvalue of $-L$ corresponding to $e_k$, 
\begin{equation}
\label{eq:LpTrace}
K_{s,r}\,
\coloneqq\,\sup_{t\in[s,r]}\sum_{k\in\N}\norm*{\Lambda_t e_k}_E\,<\,\infty. 
\end{equation}
and for $y\in E\subset H$, written $y=\sum_{k\in\N}y^ke_k$ in $H$, the following converges in $E$: 
\begin{equation}
\label{eq:LambdaLy}
L\Lambda_ty\,=\,\Lambda_tLy\,\coloneqq\,\sum_{k\in\N}y^k\Lambda_tLe_k. 
\end{equation}
\end{enumerate}
\end{assumption}

Assumption \ref{assn:Noise}~is needed for our strong It{\^o}~formula. 
We remark that the assumption is weaker than requiring $X_t\in D(L)$ or that $E=H$, so that the It{\^o}~formulas of \cite{A17,B08,I84}~do not apply to \eqref{eq:SDE}. 
See Remark \ref{remark:EX}~below. 
The condition \eqref{eq:LpTrace}~is only needed when the noise in \eqref{eq:SDE}~is white in space, but not when the noise is trace class. 

Note that with additive noise, solutions of \eqref{eq:SDE}~may be pushed out of $B(\Gamma)$ in finite time.
We therefore define the exit time of the solution from $B(\Gamma)$ as 
\begin{equation} \label{eq: tau definition}
\tau\,\coloneqq\,\min\left\{\tau_{blowup},\,\inf\left\{t>0\,:\,X_t\in E/ B(\Gamma)\right\}\right\}. 
\end{equation}
Since the geometry of $B(\Gamma)$ is generally unknown and may be very complicated, we will sometimes restrict our attention to $\Gamma_\delta$ for $\delta>0$. 
We define the exit time from $\Gamma_\delta$ as 
\begin{equation}
\label{eq:tau_delta}
\tau_\delta\,\coloneqq\,\min\left\{\tau_{blowup},\,\inf\left\{t>0\,:\,\norm*{X_t-\gamma_{\pi(X_t)}}_E=\delta\right\}\right\}. 
\end{equation}
When $\Gamma$ consists of fixed points of $\phi_t$, estimates on the distribution of  $\tau_\delta$ can be found in \cite{M23}.  
See also \cite{BG13}.  
Further work is needed to estimate the distribution of $\tau$ for more general stable invariant manifolds. 
We let $\pi_t=\pi(X_t)$ be the isochronal phase of $X_t$ for $t<\tau$.

\begin{remark}
\label{remark:EX} 
Two classes of examples of \eqref{eq:ODE}~motivate the present discussion. 
The first is the class of reaction-diffusion systems of the form 
\begin{equation}
\label{eq:RDE}
\partial_tx\,=\,(\Delta-a) x  + N(x), 
\end{equation}
where $a>0$ and $N$ is a vector of polynomials. 
With $L=\Delta-a$ and $O$ bounded, we can show that \eqref{eq:RDE}~satisfies Assumptions \ref{assn:PDE}~\&~\ref{assn:Noise}. 
For instance if $O=[0,\ell]$ with 
periodic boundaries, take $H=L^2(O)$ and $E=C(O)$, the space of continuous functions on $O$ with the supremum norm.  
In this case, it is more convenient to index the eigenpairs by $\Z$ rather than $\N$. 
We then have $e_k(x)=\sin\left(\frac{2k\pi}{\ell}x\right)$ and $\lambda_k=\frac{4\pi^2k^2}{\ell^2}$ for $k<0$, and $e_k(x)=\cos\left(\frac{1k\pi}{\ell}x\right)$ and $\lambda_k=\frac{4\pi^2k^2}{\ell^2}$ for $k\ge0$. 
Hence $\norm*{e_k}_E\le1$, and 
\[
\norm*{\Lambda_t e_k}_E\,\le\,\exp\left(-\frac{4\pi^2k^2t}{\ell^2}-at\right)\quad\text{ for $k\in\Z$ and $t>0$}. 
\]
For fixed $r>s>0$ we then have $K_{s,r}\,=\,\sup_{t\in[s,r]}\sum_{k\in\N}\norm*{\Lambda_te_k}_E\,
<\,\infty$, verifying \eqref{eq:LpTrace}. 
Assumption \ref{assn:Noise}(b) holds by the regularity of solutions to the heat equation, which can be proven as in \cite[Theorem 2.3.1]{E22}. 
The regularity of the heat equation also immediately implies \eqref{eq:LambdaLy}, since $L\Lambda_ty$ is clearly equal to $\Lambda_tLy$, as defined in \eqref{eq:LambdaLy}~for $y\in E$, and in this case the graph norm of $L$ is just $\norm*{\cdot}_{C^2}$. 
Finally, when \eqref{eq:RDE}~is perturbed by additive white noise, Assumption \ref{assn:Noise}(a)~is verified using the local inversion theorem, as in \cite{DPZ14}, since in this case the stochastic convolution $\int_0^t\Lambda_{t-s}B\,dW_s$ is contained in $E$. 

The second class of examples in which we are interested in is a class of integro-differential equations, with dynamics at $\xi\in O$ described as 
\begin{equation}
\label{eq:NFE}
\begin{aligned}
\partial_tx(t,\xi)\,&=\, -x(t,\xi) +  \int_O \omega(\xi,\zeta)f(x(t,\zeta))\,d\zeta ,  \\
\partial_ty(t,\xi)\,&=\, -\varepsilon^{-1}y(t,\xi) + x(t,\xi). 
\end{aligned}
\end{equation}
where $\omega$ and $f$ are typically bounded and Lipschitz. 
These equations are frequently employed in neuroscience (where they are typically referred to as \emph{neural field equations})~\cite{B12}~and ecology \cite{patterson2020probabilistic}. 
Mild solutions to stochastic perturbations of \eqref{eq:NFE}~only make sense when the perturbing noise is trace class \cite{KS17}, and in this case Assumption \ref{assn:Noise}~is satisfied. 
With $L=[-1,-\epsilon^{-1}]$, Assumptions \ref{assn:PDE}(a), (b), are satisfied, though (c) is not. 
However, Assumption \ref{assn:PDE}~(c) is only needed to prove the It{\^o}~formula for the isochronal phase when $W$ is not trace class, as can be seen in the arguments of Section \ref{sec:Regularity}~below. 
As stochastic neural field equations only make sense for trace class noise, the fact that \eqref{eq:NFE}~does not satisfy Assumption \ref{assn:PDE}(c) is not an issue. 

For reaction-diffusion or neural field equations, the main examples of $\Gamma$ that we have in mind are stationary or travelling waves 
\cite{AK15,kilpatrick2013wandering,maclaurin2020wandering}, and spiral waves in excitable media, as studied in \cite{BT04,BT07}~on unbounded domains and in \cite{X93}~on bounded spatial domains -- see also \cite{KMZ22}. 
We are also interested in the case of stationary patterns that remain stationary under translation or rotation 
\cite{CGG01,maclaurin2020wandering,VRFC17}. 
So long as \eqref{eq:RDE}~or \eqref{eq:NFE}~possesses a sufficiently regular invariant manifold $\Gamma$, and the parameterization $\alpha\mapsto\gamma_\alpha$ is non-degenerate, these examples satisfy Assumption \ref{assn:Manifold}. 
For instance, if $\Gamma$ consists of travelling wave solutions on a periodic domain, this is the case. 
\end{remark}

\section{Regularity of the Isochron Map \&~a Strong It{\^o}~Formula}
\label{sec:Regularity}
\begin{theorem}
\label{theorem:C2} 
Under Assumptions \ref{assn:PDE}~\&~\ref{assn:Manifold}, there exists a unique function $\pi:B(\Gamma)\rightarrow\mathcal{S}\subset\R^m$ satisfying \eqref{eq:pidefinition}. 
This function is twice continuously Fr{\'e}chet differentiable at each $x_0\in B(\Gamma)\subset E$ in the topology of $E$. 
\begin{proof}
We begin by proving that $\pi:B(\Gamma)\rightarrow\mathcal{S}$ is well-defined. 
Since $\Gamma$ is a stable manifold of \eqref{eq:ODE}, taking any sequence of positive numbers $\{\epsilon_n\}_{n\in\N}$ monotonically decreasing to zero, there exists a sequence of times $\{t_n\}_{n\in\N}$ monotonically increasing to infinity and a sequence of closed sets $U_n\subset\Gamma$ such that for all $n\in\N$ and $\gamma\in U_n$, we have 
\begin{equation}
\label{eq:phitn}
\norm*{\phi_{t_n}(x)-\gamma}_E\,\le\,\epsilon_n/2M_1, 
\end{equation}
for some constant $M_1$. 
The collection $\{\phi_{t_n}^{-1}U_n\}_{n\in\N}$ constitutes a sequence of non-empty, closed, nested subsets of the complete metric space $E$. 
Moreover, using the bi-Lipschitz condition of $\phi_t$ on $\Gamma$ and \eqref{eq:phitn}, we see that $\diam\left(\phi_{t_n}^{-1}U_n\right)\,\le\,\epsilon_n$. 
In particular, the diameter of $\phi_{t_n}^{-1}U_n$ tends to zero as $n\rightarrow\infty$. 
By Cantor's Intersection Theorem, 
there is a unique $\gamma_*$ such that $\bigcap_{n\in\N}\phi_{t_n}^{-1}U_n\,=\,\{\gamma_*\}$. 
Since 
$\alpha\mapsto\gamma_\alpha$ is invertible, 
there is a unique $\pi(x)\in\mathcal{S}$ such that $\gamma_*=\gamma_{\pi(x)}$. 

To prove that $\pi$ is $C^2$, note that our assumptions on $N$ implies that $x\mapsto\phi_t(x)$ is $C^3$ on $B(\Gamma)\cap E$ in the topology of $E$. 
Since $C^3$ is a Banach space, the map 
\[
x\,\mapsto\,\gamma_{\pi(x)}\,=\,\lim_{t\rightarrow\infty}\phi_t(x)
\]
is $C^3$ on $B(\Gamma)\cap E$. 
By assumption, $(r\mapsto\gamma_r):\mathcal{S}\rightarrow\Gamma$ is invertible with $C^2$ inverse. 
Denoting this inverse by $\gamma^{-1}$, the map $x\,\mapsto\,\pi(x)\,=\,\gamma^{-1}(\gamma_{\pi(x)})$ is $C^2$ on $B(\Gamma)\cap E$. 
\end{proof}
\end{theorem}

We 
denote the first and second derivatives of the isochron map at $x_0\in B(\Gamma)\cap E$ in the directions $x,y\in E$ by $D\pi(x_0)x$ and $D^2\pi(x_0)[x,y]$. 
The following theorem yields additional regularity of the isochron map when Assumption \ref{assn:PDE}~holds, showing in particular that its second derivative is trace class, in some sense. 

\begin{theorem}
\label{theorem:pi2Frechet}
Under Assumptions \ref{assn:PDE}~\&~\ref{assn:Manifold}, there exists $M_\pi\in(0,\infty)$ such that 
\begin{equation}
\label{eq:pi2Trace}
\begin{aligned}
\sum_{k\in\N}\norm*{D^2\pi(x_0)[e_k,e_k]}_{\R^m}\,<\,M_\pi, \qquad &\sum_{k\in\N}\norm*{D\pi(x_0)e_k}_{\R^m}\,<\,M_\pi, \\ 
\text{ and }\qquad&\sum_{k\in\N}\norm*{D\pi(x_0)e_k}_{\R^m}^2\,<\,M_\pi,
\end{aligned}
\end{equation}
uniformly in $x_0\in\Gamma_\delta$. 
Moreover, 
these quantities are Lipschitz continuous in $x_0\in\Gamma_\delta$. 
\end{theorem}

To prove Theorem \ref{theorem:pi2Frechet}, define $\tilde{\pi}:\Gamma_\delta\rightarrow\Gamma$ as $\tilde{\pi}(x)\coloneqq\gamma_{\pi(x)}$, and note that 
\begin{equation}
\label{eq:IsochronChi}
\phi_t(\tilde{\pi}(x))\,=\,\tilde{\pi}(\phi_t(x))\qquad\forall x\in\Gamma_\delta,\,\,\,t>0. 
\end{equation}
We show that \eqref{eq:pi2Trace}~is satisfied by $\pi$ replaced with $\tilde{\pi}$. 
Then, since the map $\alpha\mapsto\gamma_\alpha$ is $C^2$ with uniformly bounded first derivative, \eqref{eq:pi2Trace}~follows. 
From \eqref{eq:IsochronChi}~and Assumption \ref{assn:Manifold}, 
\begin{align}
D\tilde{\pi}(x)y\,&=\,[D\phi_t(\tilde{\pi}(x))]^{-1}D\tilde{\pi}(\phi_t(x))D\phi_t(x)y, %
\qquad\qquad\qquad\qquad \text{ and } \\ 
D^2\tilde{\pi}(x)[y,z]\,&=\,[D\phi_t(\tilde{\pi}(x))]^{-1}
\bigg(
D^2\tilde{\pi}(\phi_t(x))[D\phi_t(x)y,D\phi_t(x)z] \\
&\qquad+ 
D\tilde{\pi}(\phi_t(x))D^2\phi_t(x)[y,z] -
D^2\phi_t(\tilde{\pi}(x))[D\tilde{\pi}(x)y,D\tilde{\pi}(x)z]
\bigg). 
\nonumber
\end{align}
We now let $\norm*{\cdot}_1\coloneqq\norm*{\cdot}_{E\rightarrow E}$ and $\norm*{\cdot}_2\coloneqq\norm*{\cdot}_{E\times E\rightarrow E}$. 
Since $D\phi_t(x)^{-1}$ and $D\tilde{\pi}(\phi_t(x))$ are bounded operators, it follows that for some $C>0$ independent of $x\in B(\Gamma)$ and $t>0$, 
\begin{equation}
\label{eq:sumDpi}
\begin{aligned}
\sum_{k\in\N}\norm{D\tilde{\pi}(x)e_k}_E\,&\le\, \norm*{D\phi_t(\tilde{\pi}(x))^{-1}D\tilde{\pi}(\phi_t(x))}_1\sum_{k\in\N}\norm*{D\phi_t(x)e_k}_{E}\\
&\le\,C\sum_{k\in\N}\norm*{D\phi_t(x)e_k}_{E}. 
\end{aligned}
\end{equation}
Similarly, since $D\phi_t(\tilde{\pi}(x))^{-1}$, $D^2\tilde{\pi}(\phi_t(x))$, $D\tilde{\pi}(\phi_t(x))$, and $D^2\phi_t(x)$ are bounded uniformly in $x\in\Gamma_\delta$, for some $C_1,C_2,C_3>0$ independent of $x\in B(\Gamma)$ and $t>0$, 
\begin{equation}
\label{eq:sumD2pi}
\begin{aligned}
&\sum_{k\in\N}\norm*{D^2\tilde{\pi}(x)[e_k,e_k]}_1\,\le\, \norm*{D\phi_t(\tilde{\pi}(x))^{-1}}_E
\bigg( 
\norm*{D^2\tilde{\pi}(\phi_t(x))}_2\sum_{k\in\N}\norm*{D\phi_t(x)e_k}_E^2 \\
&\qquad\qquad\qquad
+ \norm*{D\tilde{\pi}_E(\phi_t(x))}_1\sum_{k\in\N}\norm*{D^2\phi_t(x)[e_k,e_k]}_E\\
&\qquad\qquad\qquad\qquad
+ \norm*{D^2\phi_t(\tilde{\pi}(x))}_2\sum_{k\in\N}\norm*{D\tilde{\pi}(x)e_k}_E^2
\bigg)\\
&\qquad\le\, C_1 \sum_{k\in\N}\norm*{D\phi_t(x)e_k}_E^2 + C_2\sum_{k\in\N}\norm*{D^2\phi_t(x)[e_k,e_k]}_E + C_3\sum_{k\in\N}\norm*{D\tilde{\pi}(x)e_k}^2_E. 
\end{aligned}
\end{equation}
Hence, Theorem \ref{theorem:pi2Frechet}~immediately follows as a corollary to the following result. 

\begin{lemma}
\label{prop:Trace} 
Under Assumptions \ref{assn:PDE}~\&~\ref{assn:Manifold}, there exists $M_\phi\in(0,\infty)$ such that for any $t>0$ and small $\delta>0$, 
the following hold uniformly in $x_0\in\Gamma_\delta$ and $s\in(0,t)$
\begin{equation}
\label{eq:Tr}
\begin{aligned}
&\sum_{k\in\N} \norm*{D\phi_s(x_0)[e_k]}_E\,<\,M_\phi, \\
&\sum_{k\in\N}\norm*{D\phi_s(x_0)[e_k]}_E^2\,<\,M_\phi,\quad\text{ and }\quad\sum_{k\in\N} \norm*{D^2\phi_s(x_0)[e_k,e_k]}_E\,<\,M_\phi
.
\end{aligned}
\end{equation}
Moreover, 
these quantities are Lipschitz continuous in $x_0\in\Gamma_\delta$. 
\begin{proof}
We consider the evolution equations of the operators defined in \eqref{eq:Derivatives}, written as 
\begin{equation}
\label{eq:4_Variational}
\begin{aligned}
&D\phi_t(x_0)[x]\,=\Lambda_t D\phi_0(x_0)[x] + \int_0^t\Lambda_{t-s}DN(\phi_s(x_0))D\phi_s(x_0)[x]\,ds, \\
&D^2\phi_t(x_0)[x,y]\,=\,\Lambda_tD^2\phi_0(x_0)[x,y] + \int_0^t\Lambda_{t-s}DN(\phi_s(x_0))D^2\phi_s(x_0)[x,y]\,ds \\
&\qquad\qquad\qquad\qquad+ \int_0^t\Lambda_{t-s}D^2N(\phi_s(x_0))\big[D\phi_s(x_0)[x], D\phi_s(x_0)[y]\big]\,ds, 
\end{aligned}
\end{equation}
where $DN$ and $D^2N$ are the first and second Fr{\'e}chet derivatives of $N$. 

Now, observe that $D\phi_0(x_0)[x]=x$. 
Applying Tonelli's theorem, 
\[
\begin{aligned}
\sum_{k\in\N}\norm*{D\phi_t(x_0)[e_k]}_E\,&\le\,\sum_{k\in\N}\norm*{\Lambda_te_k}_E +\sum_{k\in\N} \int_0^t\norm*{\Lambda_{t-s}DN(\phi_s(x_0))}_1\norm*{D\phi_s(x_0)[e_k]}_E\,ds\\
&\le\, \sum_{k\in\N}\norm*{\Lambda_te_k}_E + C_0\int_0^t\sum_{k\in\N}\norm*{D\phi_s(x_0)[e_k]}_E\,ds 
\end{aligned}
\]
for some $C_0>0$ which is independent of $x_0$ (using Assumption \ref{assn:PDE}, and the fact that $N$ is uniformly Lipschitz on $\Gamma_\delta$ in the topology of $E$). 
Taking $s\in(0,t)$ and letting $r\mapsto K_{s,r}$ be the non-decreasing function defined in \eqref{eq:LpTrace}, we apply Gr{\"o}nwall's inequality to obtain, for arbitrary $r>t$, 
\begin{equation}
\label{eq:Use}
\sum_{k\in\N}\norm*{D\phi_t(x_0)[e_k]}_E\,\le\, K_{s,r} e^{C_0t} \,<\,\infty. 
\end{equation} 
Furthermore, fixing $s<t<r$, \eqref{eq:Use}~implies a local integrability of the sum, 
\[
\sum_{k\in\N}\int_0^t\norm*{D\phi_s(x_0)[e_k]}_E\,ds\,=\,\int_0^t\sum_{k\in\N}\norm*{D\phi_s(x_0)[e_k]}_E\,ds\,<\,\frac{K_{s,r}}{C_0}e^{C_0t}. 
\]
We remark that the constants appearing in this bound are independent of $x_0\in\Gamma_\delta$, implying the first uniform bound \eqref{eq:Tr}. 
%

%
Similarly, observe that 
\[
\begin{aligned}
\sum_{k\in\N}\norm*{D\phi_t(x_0)[e_k]}_E^2\,&\le\,\sum_{k\in\N}\norm*{\Lambda_te_k}_E^2 + 2C_0\sum_{k\in\N}\left(\norm*{\Lambda_te_k}_E\int_0^t\norm*{D\phi_s(x_0)[e_k]}_E\,ds\right) \\
&\,\,\,+ C_0^2\sum_{k\in\N}\left(\int_0^t\norm*{D\phi_s(x_0)[e_k]}_E\,ds\right)^2. 
\end{aligned}
\]
The first and second sums are bounded, respectively, by \eqref{eq:LpTrace}~\&~\eqref{eq:Use}. 
Hence, 
\[
\begin{aligned}
\sum_{k\in\N}\norm*{D\phi_t(x_0)[e_k]}_E^2\,&\le\, K_{s,r}^2 + 2K_{s,r}^2
e^{C_0t} + C_0^2\sum_{k\in\N}\left(\int_0^t\norm*{D\phi_s(x_0)[e_k]}_E\,ds\right)^2 \\
&\le\, K_{s,r}^2\left(1+2e^{C_0t}\right) + C_0^2\int_0^t\sum_{k\in\N}\norm*{D\phi_s(x_0)[e_k]}_E^2\,ds, 
\end{aligned}
\] 
for arbitrary $s<t<r$. 
Applying Gr{\"o}nwall's inequality, we find that $\sum_{k\in\N}\norm*{D\phi_t(x_0)[e_k]}_E^2$ is finite and locally integrable in $t>0$. 
We remark that the constants appearing in the second bound of \eqref{eq:Tr}~obtained by these arguments are independent of $x_0$. 

Once we have the first two uniform bounds in \eqref{eq:Tr}, we can similarly obtain the third. 
Since $D^2\phi_0(x_0)[x,x]=x$, we have 
\[
\begin{aligned}
\sum_{k\in\N}\norm*{D^2\phi_t(x_0)[e_k,e_k]}_E\,&\le\,\sum_{k\in\N}\norm*{\Lambda_te_k}_E + C_1\int_0^t\sum_{k\in\N}\norm*{D^2\phi_s(x_0)[e_k,e_k]}_E\,ds\\
&\qquad\qquad\qquad\qquad+ C_1\int_0^t\sum_{k\in\N}\norm*{D\phi_s(x_0)[e_k]}_E^2\,ds. 
\end{aligned}
\]
By our previous results, the first and third sum in the above are finite. 
Therefore for arbitrary $s<t<r$ we have 
\[
\sum_{k\in\N}\norm*{D^2\phi_t(x_0)[e_k,e_k]}_E\,\le\, L_{s,r} + C\int_0^t\sum_{k\in\N}\norm*{D^2\phi_s(x_0)[e_k,e_k]}_E\,ds 
\]
for another constant $L_{s,r}>0$. 
We may therefore apply Gr{\"o}nwall's inequality again, to obtain the third bound in \eqref{eq:Tr}~with constants independent of $x_0$. 

To see Lipschitz continuity of the maps defined by taking $x_0$ to the values in \eqref{eq:Tr}, we remark that the first part of this proposition implies that the finite sums 
\[
\begin{aligned}
&\sum_{k=1}^K \norm*{D\phi_t(x_0)[e_k]}_E\,<\,M_\phi, \qquad\qquad\quad\,\,\,\sum_{k=1}^K \lambda_k\norm*{D\phi_t(x_0)[e_k]}_E\,<\,M_\phi, \\
&\sum_{k=1}^K\norm*{D\phi_t(x_0)[e_k]}_E^2\,<\,M_\phi, \quad
\text{ and }\qquad  \sum_{k=1}^K \norm*{D^2\phi_t(x_0)[e_k,e_k]}_E\,<\,M_\phi, 
\end{aligned}
\]
converge uniformly in $x_0\in\Gamma_\delta$ to the series in \eqref{eq:Tr}~as $K\rightarrow\infty$. 
Hence we need only prove that the finite sums are Lipschitz. 
But, from the fact that $x_0\mapsto\phi_t(x_0)$ is $C^3$ in the topology of $E$ (since $N$ is $C^4$), we know that each of the summands is Lipschitz in $\Gamma_\delta$. 
\end{proof}
\end{lemma}

We now prove a strong It{\^o}~formula for $\pi(X_t)$ that holds for $t<\tau$. 
We need the following two lemmas. 

\begin{lemma}
Let Assumptions \ref{assn:Manifold}~\&~\ref{assn:Noise}~hold. 
For each $t>0$, the It{\^o}~integral 
\begin{equation}
\label{eq:piWint}
\int_0^{t\wedge\tau}D\pi(X_s)B\,dW_s
\end{equation}
is well-defined as an $\R^m$-valued random variable. 
\begin{proof}
For a collection of independent identically distributed Brownian motions $\{\beta^k\}_{k\in\N}$, write $W_t=\sum_{k\in\N}e_k \beta^k_t$. 
Of course, this sum does not converge in any sense as an $H$-valued random variable, but as an $H_0$-valued random variable for some Hilbert space $H_0\supset H$ (see \cite[Chapter 4]{DPZ14}~for details). 
However, we need not specify $H_0$ here.  
Indeed, the integral \eqref{eq:piWint}~can be made sense of as the limit in expectation of the finite sums
\[
\sum_{k=1}^K\int_0^{t\wedge\tau}D\pi(X_s)Be_k\,d\beta_s^k\,\eqqcolon\,\sum_{k=1}^K \tilde{\beta}_t^k, 
\]
which converge in mean square. 
Using Theorem \ref{theorem:pi2Frechet}~and It{\^o}'s isometry, we obtain 
\begin{equation}
\label{eq:piWint_norm}
\begin{aligned}
\E{\norm*{\sum_{k\in\N}\tilde{\beta}_t^k}_{\R^m}}^2\,
&\le\, \sum_{k\in\N}\E{\int_0^{t\wedge\tau}\norm*{D\pi(X_s)Be_k}_{\R^m}^2\,ds}\\
&\le\, t\sup_{x\in\Gamma_\delta}\sum_{k\in\N}\norm*{D\pi(x)e_k}_{\R^m}^2\,<\,tM_\pi M_B. 
\end{aligned}
\end{equation}
\end{proof}
\end{lemma}

\begin{lemma}
\label{lemma:DpiDV}
Let Assumptions \ref{assn:Manifold}~\&~\ref{assn:Noise}~hold. 
Then, 
$D\pi(x)Ly\,\coloneqq\,\sum_{k\in\N}y^kD\pi(x)Le_k$ 
exists in $\R^m$ for all $x,y\in E$, where $\{y^k\}_{k\in\N}$ are the $\{e_k\}_{k\in\N}$ basis coefficients of $y$ in $H$. 
Moreover, $\partial_t\pi(\phi_t(x))=D\pi(\phi_t(x))V(\phi_t(x))$. 
\begin{proof}
Write $y=\sum_{k\in\N}y^ke_k$, converging in $H$. 
Note that the sum $Ly=\sum_{k\in\N}y^k\lambda_ke_k$ does not necessarily converge in $H$, but in a larger Hilbert space $H_L$.  
Recall that for $t>0$ and $x\in E$,$D\pi(x)=M(x)D\phi_t(x)$for a bounded linear operator $M(x)$ on $E$. 
Hence, $\norm*{D\pi(x)Ly}_E\le c\norm*{D\phi_t(x)Ly}_E$ for some $c>0$, and 
\[
\norm*{D\phi_t(x)Ly}_E\,\le\,\norm*{\Lambda_tLy}_E + \int_0^t\norm*{\Lambda_{t-s}DN(\phi_s(x))}_1\norm*{D\phi_s(x)Ly}_E\,ds. 
\]
By Assumption \ref{assn:Noise}(c), an application of Gr{\"o}nwell's inequality completes the proof of the convergence of $D\pi(x)Ly$. 
The statement $\partial_t\pi(\phi_t(x))=D\pi(\phi_t(x))V(\phi_t(x))$ then follows from the chain rule. 
\end{proof}
\end{lemma}

\begin{theorem}
\label{theorem:Ito} 
Let Assumptions \ref{assn:Manifold}~\&~\ref{assn:Noise}~hold. 
Then, 
\begin{equation}
\label{eq:piIto}
\begin{aligned}
\pi(X_{t\wedge\tau})\,&=\, \pi(X_0) + \int_0^{t\wedge\tau} D\pi(X_s)V(X_s)\,ds + \frac{\sigma^2}{2}\int_0^{t\wedge\tau} \sum_{k\in\N} D^2\pi(X_s)[Be_k, Be_k]\,ds \\
&\qquad + \sigma\int_0^{t\wedge\tau} D\pi(X_s)B\,dW_s,   
\end{aligned} 
\end{equation}
and all of the terms in the above expression are well-defined in $\R^m$. 
\begin{proof}
Take $t>0$ and a partition of $[0,t]$ with boundary points $\{t_k\}_{k=1}^M$. 
Let $t^*\coloneqq t\wedge\tau$ and $t_i^*\coloneqq t_i\wedge\tau$ for $i\in\{1,\ldots,M\}$. 
By Theorem \ref{theorem:C2}, $\pi(X_{t^*})$ may be written as 
\begin{equation}
\label{eq:PartitionTaylor} 
\pi(X_{t^*})\,=\,\pi(X_0) + \sum_{i=1}^M\Big( D\pi(X_{t^*_i})[X_{t^*_{i+1}}-X_{t^*_i}] + D^2\pi(w_i)[X_{t^*_{i+1}}-X_{t^*_i},X_{t^*_{i+1}}-X_{t^*_i}] \Big), 
\end{equation}
where $w_i = a_i X_{t^*_i} + (1-a_i)X_{t^*_{i+1}}$ for some $a_i\in[0,1]$. 
Since $(X_t)_{t\ge0}$ is a mild solution of \eqref{eq:SDE}, 
\[
\begin{aligned}
X_{t^*_{i+1}}-X_{t^*_i}\,&=\,\left[\Lambda_{t^*_{i+1}-t^*_i}-I\right]X_0 + \int_{t^*_i}^{t^*_{i+1}} \Lambda_{t^*_{i+1}-s}N(X_s)\,ds + \sigma\int_{t^*_i}^{t^*_{i+1}} \Lambda_{t^*_{i+1}-s}B\,dW_s \\ 
&\eqqcolon\, U^1_i + U^2_i + \sigma U_i^3. 
\end{aligned}
\]
Inserting this into \eqref{eq:PartitionTaylor}, we then have 
\begin{equation}
\label{eq:PartitionTaylorDecomp}
\begin{aligned}
\pi(X_{t^*})-\pi(X_0) \,&=\, \sum_{i=1}^M D\pi(X_{t^*_i})\left[ U^1_i \right] + \sum_{i=1}^M D\pi(X_{t^*_i})\left[ U^2_i \right] + \sigma\sum_{i=1}^M D\pi(X_{t^*_i})\left[ U_i^3\right] \\
&\quad + \sum_{i=1}^M D^2\pi(w_i)\left[U^1_i+U^2_i,U^1_i+U^2_i\right] +2\sigma\sum_{i=1}^M D^2\pi(w_i)\left[U^1_i+U^2_i,U_i^3\right] \\
&\qquad +\sigma^2\sum_{i=1}^M D^2\pi(w_i)\left[U^3_i,U^3_i\right] \\
&\eqqcolon I + II + III + IV + V + VI. 
\end{aligned}
\end{equation}

For a fixed partition $\{t_i\}_{i=1}^M$ with mesh size $h>0$, choose arbitrary $i\in\{1,\ldots,M\}$ and let $(\hat{X}_s)_{s\in[t_i,t_{i+1}]}$ be defined by the deterministic flow started at $X_{t_i}$, 
\[
\hat{X}_{t_i}\,=\,X_{t_i},\quad \hat{X}_s\,=\,\Lambda_{s-t_i}X_{t_i}+\int_{t_i}^{s}\Lambda_{s-r}N(\hat{X}_r)\,dr\quad\text{ for }\quad s\in[t_i,t_{i+1}]. 
\]
By Assumption \ref{assn:Noise}(b), note that $(\hat{X}_s)_{s\in[t_i^*,t_{i+1}^*]}$ exists as a continuous path in $D(L)$ (in the topology of the graph norm of $L$), and hence 
\[
\sup_{s\in[t_i^*,t_{i+1}^*]}\lim_{h\rightarrow0}h^{-1}\norm*{\left(\Lambda_h-I\right)\hat{X}_s}_E\,<\,\infty. 
\]
Consequently, we have the following fact, which will prove to be useful below: 
\begin{equation}
\label{eq:AssumptionF}
\sup_{s\in[t_i^*,t_{i+1}^*]}\lim_{h\rightarrow0}h^{-1/2}\norm*{\left(\Lambda_h-I\right)\hat{X}_s}_E\,=\,0. 
\end{equation}

We now show that $(I+II)\rightarrow\,\int_0^{t^*}D\pi(X_s)V(X_s)\,ds$ as $h\rightarrow0$. 
To do so, we apply the second order Taylor's theorem to $\pi(\hat{X}_{t^*_{i+1}})$ centered at $\hat{X}_{t_i}$ to obtain 
\[
D\pi(\hat{X}_{t^*_i})\left[\hat{X}_{t^*_{i+1}}-\hat{X}_{t^*_i}\right] + D^2\pi(U^5_i)\left[\hat{X}_{t^*_{i+1}}-\hat{X}_{t^*_i},\hat{X}_{t^*_{i+1}}-\hat{X}_{t^*_i}\right]\,=\,\pi(\hat{X}_{t^*_{i+1}})-\pi(\hat{X}_{t^*_i}),  
\]
where $U^5_i$ is defined for some $a_i^5\in[0,1]$ as 
\[
U_i^5\,\coloneqq\,a_i^5\hat{X}_{t^*_i} + (1-a_i^5)U^4_i. 
\]
By Gr{\"o}nwall's inequality and the Lipschitz property of $N$ (on $\Gamma_\delta$), it can be seen that 
\begin{equation}
\label{eq:diffX}
\norm*{X_s-\hat{X}_s}_E\,\sim\,O(h^2)\quad\text{ for }\,\,\,s\in[t_i,t_{i+1}]. 
\end{equation}
It then follows that 
\[
\begin{aligned}
&\norm*{\hat{X}_{t^*_{i+1}}-\hat{X}_{t^*_i} - \left[(\Lambda_{t^*_{i+1}-t^*_i}-I)X_{t^*_i} + \int_{t_i}^{t^*_{i+1}}\Lambda_{t^*_{i+1}-s}N(X_s)\,ds \right] }_E\\
&\qquad\qquad\le\, K_0h\sup_{s\in[t^*_i,t^*_{i+1}]}\norm*{\hat{X}_s-X_s}_E\,\le\,K_1h^2 
\end{aligned} 
\]
for some constants $K_0,K_1>0$. 
Using \eqref{eq:AssumptionF}, for all $t\in(t_i,t_{i+1})$ we have, 
\begin{equation}
\label{eq:hatX}
\norm*{\hat{X}_{t^*}-\hat{X}_{t^*_i}}_E\,\sim\,o(h^{1/2}). 
\end{equation}
Hence, using Lemma \ref{lemma:DpiDV}~each summand in $(I+II)$ is estimated as 
\[
\begin{aligned}
&D\pi(X_{t^*_i})\left[(\Lambda_{t^*_{i+1}-t^*_i}-I)X_{t^*_i} + \int_{t_i}^{t_{i+1}}\Lambda_{t_{i+1}-s}N(X_s)\,ds\right]\\
&\qquad=\, -D^2\pi(U^5_i)\left[\hat{X}_{t^*_{i+1}}-\hat{X}_{t^*_i},\hat{X}_{t^*_{i+1}}-\hat{X}_{t^*_i}\right] + O(h^2) + \pi(\hat{X}_{t^*_{i+1}})-\pi(\hat{X}_{t^*_i})\\
&\qquad=\,o(h) + \pi(\hat{X}_{t^*_{i+1}})-\pi(\hat{X}_{t^*_i})
\,=\,o(h) + h D\pi(\hat{X}_{t^*_i})V(\hat{X}_{t^*_i}). 
\end{aligned}
\]
Noting that $\hat{X}_{t_i^*}=X_{t_i^*}$, summing over $\{t_i^*\}_{i=1}^M$, and taking $h\rightarrow0$ yields the result. 

We now show that $III\rightarrow\,\int_0^{t^*}D\pi(X_s)B\,dW_s$. 
To do so, note that a Taylor expansion of $\Lambda_{t_{i+1}-s}$ about $s=t_{i+1}$ yields 
\begin{equation}
\label{eq:LambdaTaylor}
\Lambda_{t_{i+1}-s}\,=\, I + L\Lambda_u(t_{i+1}-s)
\end{equation}
for some $u\in(s,t_{i+1})$. 
Hence, for arbitrary $x\in\Gamma_\delta$ 
\[
\begin{aligned}
D\pi(x)\int_{t^*_i}^{t^*_{i+1}}\Lambda_{t^*_{i+1}-s} B\,dW_s \,&=\, \int_{t_i^*}^{t_{i+1}^*}D\pi(x)\left(I+(t_{i+1}-s)L\Lambda_u\right) B\,dW_s\\
&=\,D\pi(x)B[W_{t_{i+1}^*}-W_{t_i^*}]  + \int_{t_i^*}^{t_{i+1}^*} (t_{i+1}-s)D\pi(x) L\Lambda_u B\,dW_s. 
\end{aligned}
\]
By Lemma \ref{lemma:DpiDV}, $D\pi(x)L\Lambda_uB:H\rightarrow\R^m$ is a bounded linear operator, 
which we may identify with an element of $H$. 
Hence, applying the Burkholder-Davis-Gundy from \cite{MR16}, 
\[
\begin{aligned}
\ExpOp_x\left[\norm*{\int_{t_i^*}^{t_{i+1}^*}(t_{i+1}-s)D\pi(x)L\Lambda_uB\,dW_s}^2_{\R^m}\right]\,&\le\,h\ExpOp_x\left[\int_{t_i^*}^{t_{i+1}^*}\norm*{D\pi(x)L\Lambda_uB}^2_{H^m}\,ds\right]\\
&\le\, h^2\norm*{D\pi(x)L}^2\norm*{B}^2 e^{-2\omega t_i^*}. 
\end{aligned}
\]
Therefore for small $h>0$ and any $i\in\N$, 
\[
D\pi(X_{t_i^*})\int_{t^*_i}^{t^*_{i+1}}\Lambda_{t^*_{i+1}-s} B\,dW_s\,=\, D\pi(X_{t_i^*}) B[W_{t^*_{i+1}}-W_{t^*_i}] + o\left(h^2\right). 
\]
Observing that $K=\floor{t/h}$, we have the following in probability
\[
\begin{aligned}
III\,&=\, \sum_{i=0}^KD\pi(X_{t^*_i})\left[\int_{t^*_i}^{t^*_{i+1}} \Lambda_{t^*_{i+1}-s}B\,dW_s\right]\\
&=\,\sum_{i=0}^K D\pi(X_{t^*_i})B[W_{t^*_{i+1}}-W_{t^*_i}] + o\left(h^{1/2}\right)\,\xrightarrow[h\rightarrow0]{}\,\int_0^{t^*}D\pi(X_s)B\,dW_s. 
\end{aligned}
\]

We now show that $VI\rightarrow\int_0^{t^*}\sum_{k\in\N}D^2\pi(X_s)\left[Be_k,Be_k\right]\,ds$ as $h\rightarrow0$. 
By \eqref{eq:LambdaTaylor}, 
\begin{equation}
\label{eq:help1}
\begin{aligned}
&\int_{t^*_i}^{t^*_{i+1}}\sum_{k\in\N}D^2\pi(w_i)\left[ \Lambda_{t^*_{i+1}-s}Be_k, \, \Lambda_{t^*_{i+1}-s}Be_k  \right]\,ds\\
&\qquad=\, \int_{t^*_i}^{t^*_{i+1}}\sum_{k\in\N}\bigg(D^2\pi(w_i)\left[Be_k,Be_k\right]  + 2D^2\pi(w_i)\left[L \Lambda_uBe_k,Be_k\right](t^*_{i+1}-s) \\
&\qquad\qquad+ D^2\pi(w_i)\left[L\Lambda_uBe_k,L \Lambda_uBe_k\right](t^*_{i+1}-s)^2\bigg)\,ds\\
&\qquad=\, (t^*_{i+1}-t^*_i)\sum_{k\in\N}D^2\pi(w_i)\left[Be_k,Be_k\right] + O(h^2). 
\end{aligned}
\end{equation}
It{\^o}'s isometry and \eqref{eq:help1}~then imply that 
\begin{equation}
\label{eq:help2}
\begin{aligned}
&\Bigg|\ExpOp\Bigg[D^2\pi(w_i)\left[\int_{t_i^*}^{t_{i+1}^*}\Lambda_{t_{i+1}^*-s}B\sum_{k\in\N}e_k\,d\beta_s^k,\,\int_{t_i^*}^{t_{i+1}^*}\Lambda_{t_{i+1}^*-s}B\sum_{\ell\in\N}e_\ell\,d\beta_s^\ell,\right]\\
&\qquad\qquad\qquad\qquad\qquad\qquad-\int_{t_i^*}^{t_{i+1}^*}\sum_{k\in\N}D^2\pi(w_i)\left[\Lambda_{t_{i+1}^*-s}Be_k,\,\Lambda_{t_{i+1}^*-s}Be_k\right]\,ds\Bigg]\Bigg|\\
&\le\,\norm*{\pi}_{C^2}\abs*{\sum_{k\in\N}\ExpOp\left[\norm*{\int_{t_i^*}^{t_{i+1}^*}\Lambda_{t_{i+1}^*-s}Be_k\,d\beta_s^k}^2-\int_{t_i^*}^{t_{i+1}^*}\norm*{\Lambda_{t_{i+1}^*-s}Be_k}^2\,ds\right]}\,=\,0, 
\end{aligned}
\end{equation}
%
where we have used \eqref{eq:help1}~to justify exchanging the sums and expectation in \eqref{eq:help2}. 
Taking \eqref{eq:help1}~\&~\eqref{eq:help2}~together, summing over $k\in\{1,\ldots K\}$, and letting $h\rightarrow0$, we have convergence in probability of $VI$. 
Finally, $IV$ and $V$ tend to zero as $h\rightarrow0$ by \eqref{eq:AssumptionF}~and \eqref{eq:diffX}. 
\end{proof}
\end{theorem}

\bibliographystyle{amsplain}
\bibliography{bibliography}

\subsubsection*{Acknowledgements}
The author would like to thank James MacLaurin, for many valuable discussions, Wilhelm Stannat, for his insightful critiques, and J{\"u}rgen Jost, for his patience and guidance. 

\end{document}